\documentclass[12pt,a4paper]{article}
\usepackage{amsmath,amssymb,amsxtra,amsthm,amsbsy,mathtools}
\usepackage{mathpazo}

\DeclarePairedDelimiter{\diagfences}{(}{)}

\newcommand{\diag}{\operatorname{diag}\diagfences}
\theoremstyle{plain}
\newtheorem{thm}{Theorem}
\newtheorem*{pf}{Proof}

\title{\Large \bf
Seemingly unrelated and fixed-effect panel regressions: collinearity and singular dispersion}

\author{Harry Haupt\footnote{Address correspondence to: Harry Haupt,
Chair of Statistics and Data Analytics, University of Passau, Germany.
harry.haupt@wiwi.uni-passau.de}}

\begin{document}
\maketitle
\thispagestyle{empty}

\begin{abstract}
The paper discusses identification conditions and representations of generalized least squares estimators of regression parameters in multivariate linear regression models such as seemingly unrelated and fixed effect panel models. Results are presented on identification for unrestricted dispersion structure, general heteroskedasticity and cross-equation dependence, considering explicit and implicit linear restrictions, singular dispersion matrix and collinear design matrix.
\end{abstract}

\begin{description}
\item[Keywords:] Least squares, linear restrictions, seemingly unrelated regression, \\
fixed effects, singular dispersion matrix, collinearity.
\end{description}

\vskip1cm

\thispagestyle{empty}

\newpage
\section{Model structures and assumptions}

Consider regression equations
$y_{t,i} = \mu_{t,i} + u_{t,i}$, $1 \le i \le n$, $1 \le t \le m$,
where $y_{t,i}$ is the random response variable, $\mu_{t,i}=E(y_{t,i})$ is the unknown systematic component and $u_{t,i}$ is the random error variable.
The equations
\begin{equation*}
\mu_{t,i}=\sum_{k=1}^{K_i} \beta_{k,i} x_{t,k,i}, \quad 1 \le i \le n, \quad 1 \le t \le m,
\tag{SUR}
\end{equation*}
where $x_{k,i}$ are known covariates and $\beta_{k,i}$ are fixed but unknown parameters, is called seemingly unrelated regressions (SUR) model. This model class, introduced by Zellner (1962), found many applications in mathematical statistics, econometrics and related disciplines and elicits an ongoing research interest (e.g., Rao, 1975, Mandy \& Martins-Filho, 1993, Jackson, 2002, Qian, 2008, Kurata \& Matsura, 2016, Zeebari et al., 2018, or Hou \& Zhao, 2019).

For every $t$, model (SUR) can be written as
\[
\mathbf{y}_{t,\bullet} = \mathbf{X}_{t,\bullet} \boldsymbol{\beta} + \mathbf{u}_{t,\bullet},
\quad 1 \le t \le m,
\]
where $\mathbf{y}_{t, \bullet}$ is the $(n \times 1)$ vector of response variables (in time $t$),
$\mathbf{X}_{t,\bullet} = \diag{\mathbf{x}'_{t,\bullet,1}, \ldots, \mathbf{x}'_{t,\bullet,n}}$ is an $(n \times K)$ matrix, $\mathbf{x}'_{t,\bullet,i}$ represents the $(1 \times K_i)$ vector of covariates in equation $i$ at time $t$, and $\boldsymbol{\beta}$ is a $(K \times 1)$ vector of parameters, where $K=\sum_iK_i$, and $\mathbf{u}_{t,\bullet}$ is an $(n \times 1)$ error vector.

Stacking the equations for every $t$ leads to the compact representation
\begin{equation}
\mathbf{y} = \mathbf{X} \boldsymbol{\beta} + \mathbf{u},
\label{eq:SUR}
\end{equation}
where $\mathbf{y}$ is the $(T \times 1)$ vector of response variables and $T=mn$, $\mathbf{X}$ is the $(T \times K)$ matrix of covariates, $\boldsymbol{\beta}$ is the $(K \times 1)$ vector of unknown parameters and $\mathbf{u}$ is the $(T \times 1)$ vector of error variables, and $T > K$.
For our general discussion of least squares theory in Section 2 we assume the following general stochastic specification of model (1):
All covariates are non-stochastic,
$E(\mathbf{u}) = \mathbf{0}$,
and the system dispersion structure is represented by the unrestricted and known (e.g., Rao, 1975, Hou \& Zhao, 2019) $(T \times T)$ matrix
$E(\mathbf{{u}} \mathbf{{u}}') \equiv \sigma^2\mathbf{\Omega}$,
where $\sigma^2$ is an unknown positive and finite scalar parameter, $\mathbf{\Omega}$ is a symmetric, nonnegative definite matrix with $tr(\mathbf{\Omega})= T$.
The resulting triplet is commonly denoted as general Gauss-Markoff model:
\[
\mathcal{A} = \{\mathbf{y},\mathbf{X}\boldsymbol{\beta}, \sigma^2 \mathbf{\Omega}\}.
\]

In model $\mathcal{A}$, both $\mathbf{X}$ and $\mathbf{\Omega}$ may be deficient in rank. As a consequence, $y \in \mathcal{M}(\mathbf{X} : \mathbf{\Omega})$ with probability one (see Rao, 1973, Lemma 2.1), where $\mathcal{M}(\mathbf{A})$ denotes the linear manifold generated by the columns of the matrix $\mathbf{A}$. Clearly, for the classical situation where $\mathbf{\Omega}$ has full rank, $y \in \mathcal{M}(\mathbf{\Omega})$ with probability one, as then the complete $\mathbb{R}^T$ is spanned (see Haupt \& Oberhofer, 2002).
The problem of a singular dispersion arises frequently in applications involving adding-up restrictions such as demand or share equation models (e.g., Barten, 1969, or Bewley, 1986). Further, from an algebraic point of view, the singular dispersion matrix problem lies at the very heart of classical regression theory due to the singular (idempotent and orthogonal) projection matrices of least squares. The collinearity problem on the other hand has an impressive track record in the statistical literature, with an even increasing intensity due to the omnipresent phenomenon of big data (and covariate) availability (see Dormann et al., 2013, for a recent review).

The aim of this paper is twofold.

First, we provide a thorough discussion of results on identification conditions and representations of least squares estimators and their relations in the general framework defined by model (1) and special cases of this model class such as panel data models.
Since the seminal work of Rao (1965, 1971, 1973) on unified least squares theory, the model $\{\mathbf{y},\mathbf{X}\boldsymbol{\beta}, \sigma^2 \mathbf{\Omega}\}$ has been subject of intensive discussion in mathematical statistics and econometrics.
Various works consider ordinary or generalized least squares (OLS or GLS) or feasible GLS estimation of $\boldsymbol{\beta}$ (or linear transformations thereof) under different constellations of collinear design matrix ($rk(\mathbf{X})=P \le K$), and/or singular dispersion matrix ($rk(\mathbf{\Omega})=M \le T$). A natural remedy to deal with collinearity is to reduce the dimensionality of the parameter space by introducing non-sample information in form of linear parametric restrictions (e.g., Rao, 1975, Hill \& Adkins, 2001, Zeebari et al., 2018). The latter, on the other hand, are a direct consequence of a singular dispersion.
Hence, it is quite natural to include such restrictions as well as content-related restrictions (e.g., symmetry, exclusion, etc.) from economics, politics, etc. (see e.g., Berndt \& Savin 1975, Bewley, 1986, or Jackson, 2002, for examples) and discuss the restricted model
\[
\mathcal{A}_r =\{\mathbf{y},\mathbf{X}\boldsymbol{\beta}|\mathbf{R}\boldsymbol{\beta}
= \mathbf{r}, \sigma^2 \mathbf{\Omega}\};
\]
e.g., Theil (1971), Rao \& Mitra (1971), Lawson \& Hanson (1974), Kreijger \& Neudecker (1977), Baksalary \& Kala (1981), Harville (1981), Magnus \& Neudecker (1988), Haupt \& Oberhofer (2002, 2006a,b) and, more recently, the works of Tian (e.g., Tian, 2007, 2009, Tian et al., 2008, Tian \& Puntanen, 2009, Tian \& Wiens, 2006, Tian \& Zhang, 2011).

Second, we apply the key results to GLS estimation in the fixed effects (FE) panel data regression framework. We provide a general result on the equality of generalized least squares estimators under assumptions (e.g., Mandy \& Martins-Filho, 1993, Jackson, 2002, or Haupt \& Oberhofer, 2006a) on the dispersion structure frequently used in applications:
The error variables are heteroskedastic both in every equation and across equations,
contemporaneously correlated and generally uncorrelated in time,
$Cov(u_{t,i},u_{s,j}) = 0$, for $t \ne s, \,\, 1 \le i,j \le n$,
$E(\mathbf{u}_{t,\bullet}\mathbf{u}_{t,\bullet}') \equiv \mathbf{\Sigma}_t$, $1 \le t \le m$,
where $\mathbf{\Sigma}_1, \ldots, \mathbf{\Sigma}_m$ are known symmetric, nonnegative definite $(n \times n)$ matrices with elements
$Cov(u_{t,i},u_{t,j}) \equiv \sigma_{t,i,j}$, $1 \le t \le m, \,\,1 \le i,j \le n$.
Then the system dispersion structure is represented by the $(T \times T)$ matrix
$E(\mathbf{{u}} \mathbf{{u}}') = \diag{\mathbf{\Sigma}_{1}, \ldots, \mathbf{\Sigma}_{m}}$.
The results presented for the one-way FE model hold for more general models (such as two-way FE models) and any block-diagonal dispersion structure, for example a Kronecker structure (e.g., Kiefer, 1980), following from the additional assumption that contemporaneous correlations do not vary over time, implying
$\mathbf{\Sigma}_t = \mathbf{\Sigma}$, $\forall t$,
and hence
$E(\mathbf{{u}} \mathbf{{u}}') = \mathbf{I}_m \otimes \mathbf{\Sigma}$.

The remainder of the paper is organized as follows: Section 2 starts with a brief discussion of the classical textbook case to introduce familiar notation and concepts under the assumptions of regular design and dispersion. Subsections 2.1-2.3 discuss rank conditions for least squares estimation of SUR models under different rank conditions and provides relationships between the resulting estimators.
Maintaining regularity of the dispersion matrix, Section 2.1 discusses the basic structure of the design matrix when addressing the problem of (near) collinearity and a general concept to treat exact and non-exact non-sample information as a remedy.
Under the assumption of a singular dispersion matrix, Section 2.2 states the crucial rank condition required for GLS estimators using a Moore-Penrose inverse.
Section 2.3 generalizes to the case of collinear design matrix and singular dispersion matrix.
Section 3 provides a general result for the within estimator for the fixed effects panel data regression model and discusses its relationship to previous results.

\section{Generalized least squares and restrictions}

The following brief discussion of the standard textbook case, where $rk(\mathbf{X})=P=K$ (regular design matrix) and $rk(\mathbf{\Omega})=M=T$ (regular dispersion matrix), is used to introduce notation and concepts from linear algebra and multivariate statistics.
Under the assumptions of model (1) stated above, the ordinary least squares (OLS) estimator of $\boldsymbol{\beta}$ is given by
\begin{equation*}
\hat{\boldsymbol{\beta}}_{OLS}
=\mathbf{(X'X)}^{-1}\mathbf{X'y}
 = \mathbf{X}^+ \mathbf{y},
\end{equation*}
where $\mathbf{X}^+ $ denotes the Moore-Penrose generalized inverse.
For convenience we restate the conditions satisfied by the Moore-Penrose inverse of $\mathbf{X}$:
$\mathbf{X} \mathbf{X}^+ \mathbf{X} = \mathbf{X}$,
$(\mathbf{X} \mathbf{X}^+)' = \mathbf{X} \mathbf{X}^+$,
$\mathbf{X}^+ \mathbf{X} \mathbf{X}^+ = \mathbf{X}^+$, and
$(\mathbf{X}^+ \mathbf{X})' = \mathbf{X}^+ \mathbf{X}$.
$\mathbf{X}^+$ is unique and $rk(\mathbf{X}^+) = rk(\mathbf{X})$. Further,
 $\mathbf{X} \mathbf{X}^+, \mathbf{I}_T - \mathbf{X} \mathbf{X}^+,
\mathbf{X}^+ \mathbf{X}$ and $\mathbf{I}_K - \mathbf{X}^+ \mathbf{X}$ are idempotent,
$(\mathbf{X}^+)^+  = \mathbf{X}$, and $(\mathbf{X}^+)' = (\mathbf{X}')^+$.
The OLS residuals are equal to
$\hat{\mathbf{u}}
= \mathbf{y} - \mathbf{X} \hat{\boldsymbol{\beta}}_{OLS}
= \mathbf{y} - \mathbf{X}\mathbf{X}^+\mathbf{y}
= (\mathbf{I}_T-\mathbf{X}\mathbf{X}^+)\mathbf{y}$,
where
$\mathbf{P_X} \equiv \mathbf{X}\mathbf{X}^+=\mathbf{X}(\mathbf{X}' \mathbf{X})^{-1} \mathbf{X}'$ and
$\mathbf{M_X} \equiv \mathbf{I}_T-\mathbf{X}\mathbf{X}^+=\mathbf{I}_T-\mathbf{X}(\mathbf{X}' \mathbf{X})^{-1} \mathbf{X}'$
are the orthogonal projection matrices of OLS.
Under the assumptions stated above, $\hat{\boldsymbol{\beta}}_{OLS}$ is well known to be linear unbiased but not efficient, as it ignores the non-spherical structure of the dispersion matrix.
The GLS (generalized least squares) estimator of $\boldsymbol{\beta}$ in (1) is best linear unbiased and is, using the definition $\mathbf{C} \equiv \mathbf{X}' \boldsymbol{\Omega}^{-1} \mathbf{X}$, given by
\begin{equation*}
\hat{\boldsymbol{\beta}}_{GLS}
=\mathbf{C}^{-1} \mathbf{X}'\boldsymbol{\Omega}^{-1} \mathbf{y}.
\end{equation*}
Although the assumption that $\mathbf{X}$ is regular holds in many applications, we frequently encounter situations where $\mathbf{X}'\mathbf{X}$ has some very small roots (see Section 2.1). To deal with such cases of near collinearity, an option is to introduce non-sample a priori information in order to reduce the dimensionality of the parameter space by imposing linear restrictions on the parameter vector $\boldsymbol{\beta}$,
\begin{equation}
\mathbf{R}\boldsymbol{\beta}=\mathbf{r}.
\label{eq: explicitrestrictions}
\end{equation}
where both the $(q\times K)$ matrix $\mathbf{R}$ and the $(q \times 1)$ vector $\mathbf{r}$ are known.
We assume that the restrictions (2) are consistent (i.e. non-redundant) and satisfy the rank condition
\begin{equation}
rk(\mathbf{R}) = rk(\mathbf{R},\mathbf{r}),
\end{equation}
as well as the identification condition
\begin{equation}
\left( \begin{array}{c} \mathbf{R} \\ \mathbf{X} \end{array} \right) \,\, \textrm{has full column rank}.
\end{equation}
The restricted OLS (ROLS) estimator can be derived by reparametrization (see Rao, 1965, 4a.9, or Haupt \& Oberhofer, 2002):
First, let the dispersion matrix be diagonalized according to $\mathbf{\Omega} = \mathbf{F}^* \boldsymbol{\Lambda}^* {\mathbf{F}^*}'$,
where the orthogonal $(T \times T)$ matrix ${\mathbf{F}}^*$ contains the eigenvectors of $\mathbf{\Omega}$ and $\boldsymbol{\Lambda}^*$ its eigenvalues. It is well known that $\hat{\boldsymbol{\beta}}_{GLS}$ can also be obtained from reparametrization by applying OLS to the model
$\mathbf{y}^* = \mathbf{X}^*\boldsymbol{\beta} + \mathbf{u}^*$,
where $\mathbf{y}^* ={\boldsymbol{\Lambda}^*}^{-1/2}{\mathbf{F}^*}'\mathbf{y}$ and $\mathbf{X}^* ={\boldsymbol{\Lambda}^*}^{-1/2}{\mathbf{F}^*}'\mathbf{X}$.

Inverting the linear restrictions (2) leads to
\[
\boldsymbol{\beta} = \mathbf{R}^+\mathbf{r} + \mathbf{N}_R\mathbf{c},
\]
where $\mathbf{R}^+$ is the Moore-Penrose inverse of $\mathbf{R}$, $\mathbf{R}^+\mathbf{r}$ can be interpreted as a particular solution to (2), and the columns of the matrix $\mathbf{N}_R$ are the basis vectors spanning the null space on $\mathbf{R}$.
Then, model (1) can be re-written according to
$E(\mathbf{y}) = \mathbf{X} \boldsymbol{\beta} = \mathbf{X}\mathbf{R}^+\mathbf{r} + \mathbf{X}\mathbf{N}_R\mathbf{c}$,
and the ROLS estimator
\begin{eqnarray}
\hat{\boldsymbol{\beta}}_{ROLS}
= \hat{\boldsymbol{\beta}}_{OLS}
+ \mathbf{(X'X)}^{-1}\mathbf{R'}[\mathbf{R}\mathbf{(X'X)}^{-1}\mathbf{R'}]^{-1} (\mathbf{R}\hat{\boldsymbol{\beta}}_{OLS} - \mathbf{r}),
\end{eqnarray}
can be derived from applying OLS to the model $\mathbf{y}^{**} = \mathbf{X}^{**} \mathbf{c} + \mathbf{u}^{**}$, where
$\mathbf{y}^{**} = \mathbf{y} - \mathbf{XR}^+\mathbf{r}$ and $\mathbf{X}^{**}=\mathbf{XN}_R$.
The restricted GLS (RGLS) estimator, given by
\begin{equation}
\hat{\boldsymbol{\beta}}_{RGLS} =
\hat{\boldsymbol{\beta}}_{GLS} + \mathbf{C}^{-1}\mathbf{R'(R} \mathbf{C}^{-1} \mathbf{R')}^{-1} (\mathbf{R}\hat{\boldsymbol{\beta}}_{GLS} - \mathbf{r}),
\end{equation}
can be derived analogously by applying the reparametrizations (labelled by * and **, respectively) described above.
Estimators derived under the regularity conditions employed in this section, are discussed in detail in texts on regression analysis (see e.g., Amemiya, 1985, for an authoritative treatment).

\subsection{(Almost) collinear design and regular dispersion}

The case of a collinear design matrix $rk(\mathbf{X})=P<K$ or highly correlated covariates is one of the most puzzling topics in modern statistics. The problem has a long history in mathematics and statistics and its discussion gains further momentum with modern techniques in machine learning and its applications to statistics. Many dimensions of the problem as well as a voluminous literature are discussed in the recent survey of Dormann et al. (2013). In the case of (near) linear dependency of covariates, a natural remedy is to introduce non-sample information in form of linear restrictions on the parameters in model (1). See Hill \& Adkins (2001) for a survey of methods and contributions and Haupt \& Oberhofer (2005) for a general treatment of non-exact non-sample information in model (1).
For every $n$, write model (SUR) as
\[
\mathbf{y}_{\bullet,i} =
\mathbf{X}_{\bullet,i} \boldsymbol{\beta}_{\bullet,i}  + \mathbf{u}_{\bullet,i},
\quad 1 \le i \le n,
\]
where $\mathbf{y}_{\bullet,i}$ is the $(m \times 1)$ vector of response variables,
$\mathbf{X}_{\bullet,i}$ is an $(m \times K_i)$ covariate matrix, $\boldsymbol{\beta}_{\bullet,i}$
contains the $K_i$ parameters, and $\mathbf{u}_{\bullet,i}$ is an $(m \times 1)$ error vector.
Stacking this system leads to $\mathbf{y} = \mathbf{X} \boldsymbol{\beta} + \mathbf{u}$, where $\mathbf{X} = \diag{\mathbf{X}_{\bullet,1}, \ldots, \mathbf{X}_{\bullet,n}}$.
Then the generalized ridge estimator can be represented by
\begin{equation*}
\hat{\boldsymbol{\beta}}_{ridge}
=(\mathbf{X'X} + \boldsymbol{\Psi})^{-1}\mathbf{X'y},
\end{equation*}
while the ordinary ridge estimator is based on the assumption of an equation-specific scalar ridge component $\psi_i$; then $\boldsymbol{\Psi} = \diag{\psi_1\mathbf{I}_{K_1}, \ldots, \psi_n\mathbf{I}_{K_n}}$ is a $(K \times K)$ diagonal matrix containing equation-specific ridge parameters  $\psi_i>0$, $\forall i$.
Diagonalizing $\mathbf{X}'\mathbf{X}$ by $\mathbf{F}_X\boldsymbol{\Lambda}_X\mathbf{F}_X'$, we can re-write the ridge regression design matrix as
$\mathbf{X'X} + \boldsymbol{\Psi}
= \mathbf{F}_X\boldsymbol{\Lambda}_X\mathbf{F}_X' + \mathbf{F}_X\boldsymbol{\Psi}\mathbf{F}_X'
= \mathbf{F}_X(\boldsymbol{\Lambda}_X + \boldsymbol{\Psi})\mathbf{F}_X'$.
Thus, even in the case where the design is collinear, all roots are bounded away from zero.
See Rao (1975) for a general treatment on $\hat{\boldsymbol{\beta}}_{ridge}$ and Zeebari et al. (2018) on its feasible estimation.
Note that $\hat{\boldsymbol{\beta}}_{ridge}$ is a special case of the estimator proposed by Haupt \& Oberhofer (2005), resulting from GLS estimation of
\begin{equation*}
\left( \begin{array}{c} \mathbf{y} \\ \mathbf{r} \end{array} \right)
=
\left( \begin{array}{c} \mathbf{X} \\ \mathbf{RX}_f \end{array} \right)
\boldsymbol{\beta}
+
\left( \begin{array}{c} \mathbf{u} \\ \mathbf{v} \end{array} \right),
\qquad
\left( \begin{array}{c} \mathbf{u} \\ \mathbf{v} \end{array} \right) \sim
\left( \left( \begin{array}{c} \mathbf{0} \\ \mathbf{0} \end{array} \right),
\left( \begin{array}{cc} \sigma^2\boldsymbol{\Omega} & \mathbf{0} \\ \mathbf{0} & \boldsymbol{\Theta} \end{array} \right)
\right),
\end{equation*}
allowing for additional incomplete and stochastic information in model (1).
It addresses a situation typically occuring in prediction or missing values problems, where additional observations of the covariates $\mathbf{X}_f$ are available, but the corresponding response observations $\mathbf{y}_f$ are unobservable.
If no additional observations are available, this estimator is equal to the mixed estimator proposed by Theil \& Goldberger (1961) and Theil (1963) by setting $\mathbf{X}_f = \mathbf{I}_q$. If in addition the restrictions are exact (i.e. $\mathbf{v}=\mathbf{0}$), the estimator is equal to
$\hat{\boldsymbol{\beta}}_{RGLS}$ in equation (6).

\subsection{Regular design and singular dispersion}

The case given by $rk(\mathbf{X})=P=K$ and $rk(\mathbf{\Omega})=M<T$ (singular dispersion matrix) arises naturally in least squares estimation and in many applications dealing with adding-up conditions (see Haupt \& Oberhofer, 2002, 2006a,b, and the literature cited therein).
In this case $\boldsymbol{\Omega}$ has $M$ positive eigenvalues and $T-M$ eigenvalues equal to zero.
For all $t= 1, \ldots, m$ consider
${\mathbf{\Sigma}_t} {\mathbf{F}_t}^* = {\mathbf{F}_t}^* {\boldsymbol{\Lambda}_t}^*$,
where the columns of the orthogonal $(n \times n)$ matrix ${\mathbf{F}_t}^*$ are constituted by the eigenvectors of ${\mathbf{\Sigma}_t}$ and $\boldsymbol{\Lambda}_t$ is the diagonal matrix containing its eigenvalues.
Define a partition ${\mathbf{F}_t}^* = ({\mathbf{A}_t}, {\mathbf{F}_t})$, where ${\mathbf{A}_t}$ is the matrix of eigenvectors corresponding to the zero eigenvalues of ${\boldsymbol{\Sigma}_t}$ and ${\mathbf{F}_t}$ is the matrix of eigenvectors corresponding to the positive eigenvalues.
Further, define ${\mathbf{\Lambda}_t}^* = \diag{\mathbf{0}, {\boldsymbol{\Lambda}_t}}$, where ${\boldsymbol{\Lambda}_t}$ is a diagonal matrix of the positive eigenvalues of $\mathbf{{\Sigma}_t}$.
Then, ${\mathbf{A}_t}$ is an orthogonal base of the null space on ${\boldsymbol{\Sigma}_t}$.
Using analogous notation, for the stacked system (1) let
\begin{equation}
\mathbf{\Omega A}=\mathbf{0},
\end{equation}
where $\boldsymbol{\Omega}$ can be diagonalized as
$\boldsymbol{\Omega} = \mathbf{F}^* \boldsymbol{\Lambda}^* {\mathbf{F}^*}'= \mathbf{F} \mathbf{\Lambda} \mathbf{F}'$
with Moore-Penrose inverse
\begin{equation}
\boldsymbol{\Omega}^+ = \mathbf{F} \mathbf{ \Lambda}^{-1} \mathbf{F}'.
\end{equation}
Then, premultiplication of (1) by the nonsingular matrices
$\mathbf{A}'$ and $\mathbf{ \Lambda}^{-1/2} \mathbf{F}'$, respectively, leads to the transformed model
\begin{eqnarray}
\mathbf{A}'\mathbf{y}
&=& \mathbf{A}' \mathbf{X}\boldsymbol{\beta} + \mathbf{A}'\mathbf{u},
\\
\boldsymbol{ \Lambda}^{-1/2} \mathbf{F}'\mathbf{y}
&=& \boldsymbol{ \Lambda}^{-1/2} \mathbf{F}'\mathbf{X}\boldsymbol{\beta}
  + \boldsymbol{ \Lambda}^{-1/2} \mathbf{F}' \mathbf{u}.
\end{eqnarray}
From (7) follows
$E (\mathbf{A}' \mathbf{uu}'\mathbf{A}) = \sigma^2 \mathbf{A}' \boldsymbol{\Omega}\mathbf{A} = \mathbf{0}$,
and, as a consequence, in (9),
\begin{equation}
\mathbf{A}' \mathbf{u} = \mathbf{0},
\qquad
\mathbf{A}' \mathbf{y} = \mathbf{A}' \mathbf{X}\boldsymbol{\beta}
\end{equation}
hold with probability one. The latter condition is equal to assuming a singular dispersion matrix,
as $\mathbf{A}'\mathbf{y}$ is equal to a vector of constants with probability one.
Haupt \& Oberhofer (2002) discuss this case (for $M \le T-1$) and consider the resulting implicit parameter restrictions
\begin{equation}
\mathbf{G} \boldsymbol{\beta} = \mathbf{g} \quad a.s.
\end{equation}
As discussed in the regular case above, from (8) and the orthogonality of $\mathbf{F}$ follows that the transformed regression model (10) has spherical disturbances.
Then the best linear unbiased estimator can be obtained by OLS estimation of the transformed model (10) subject to the linear restrictions on $\boldsymbol{\beta}$ in (12), which are due to the singularity of the dispersion matrix.
As $\mathbf{F}$ and $\mathbf{\Lambda}$ are nonsingular, this estimator is feasible whenever the design matrix is regular, i.e. $P=K$ (e.g., Rao, 1965, and Theil, 1971).

An alternative estimation strategy in this case is to employ the Moore-Penrose inverse of $\mathbf{\Omega}$.
The Moore-Penrose inverse least squares estimator is given by
\begin{equation}
\hat{\boldsymbol{\beta}}_{MLS} = \mathbf{C}_+^{-1} \mathbf{X'\Omega^+ y}.
\end{equation}
The true inverse
\begin{equation}
\mathbf{C}_+^{-1} = (\mathbf{X}' \boldsymbol{\Omega}^+ \mathbf{X})^{-1}
= \mathbf{(X' \mathbf{F}} \boldsymbol{\Lambda}^{-1} \mathbf{F}'\mathbf{X})^{-1}.
\end{equation}
exists if the matrix $\mathbf{F}'\mathbf{X}$ has full column rank, which is known as Theil's (1971, Theorem 6.6, Assumption 6.3) first rank condition.

The following result (Haupt \& Oberhofer, 2000) provides necessary and sufficient conditions for the rank condition of Theil and improve on Theil's original suggestion and the subsequent work of Dhrymes \& Schwarz (1987). For ease of exposition we consider the common case in statistical applications, where the eigenvector $\mathbf{a}$ corresponds to the single zero root of $\boldsymbol{\Sigma}_t$ for every $t$ .

\begin{thm}Theil's rank condition is violated, if there are collinearities in any of
the $n$ equations, or, if there exists a nonzero linear combination of the
covariates of the $i$th equation, for all equations $i$ with nonzero weight $a_i$ in
$\sum_i {a}_i\mathbf{X}_{\bullet,i}\boldsymbol{\beta}_{\bullet,i} = \mathbf{s}$
(or, equivalently,
$\mathbf{a}'\mathbf{X}_{t,\bullet}\boldsymbol{\beta} =s_t$).
\end{thm}

\begin{pf}
Equivalent to a column rank deficit of the matrix $\mathbf{F}'\mathbf{X}$ is the existence of a nonzero vector $\mathbf{d}$ fulfilling
\begin{equation}
\mathbf{F}'_{t} \mathbf{X}_{t,\bullet} \mathbf{d = 0}, \qquad \forall t.
\end{equation}
Due to the orthogonality of $\mathbf{a}$ and $\mathbf{F}_{t}$, equation (15) implies
$\mathbf{X}_{t,\bullet} \mathbf{d}
= \mathbf{a}{\mathbf{a}}' \mathbf{X}_{t,\bullet} \mathbf{d}$, $\forall t$,
i.e.\ the existence of $T$ scalars $s_t$ with
\begin{equation}
\mathbf{X}_{t,\bullet} \mathbf{d} = \mathbf{a}s_t, \qquad \forall t.
\end{equation}
By using the notation of Subsection 2.1 and conducting row manipulations, equation (16) can be rendered to
\begin{equation}
\mathbf{X}_{\bullet,i} \mathbf{d}_i = a_i \mathbf{s}, \qquad \forall i,
\end{equation}
where $\mathbf{d}_i$ contains a suitable selection of the elements in $\mathbf{d}$.
Two cases arise:
First, let us assume the existence of an equation $j, (1 \le j \le n)$ with
$rk(\mathbf{X}_{\bullet,j}) < K_j$. Then
there exists a nonzero vector $\mathbf{f}$ with $\mathbf{X}_{\bullet,j}\mathbf{f = 0}$.
Now let $\mathbf{d}_i = \mathbf{0}$ for $i \ne j$, $\mathbf{d}_j = \mathbf{f}$.
Hence we have found a nonzero vector $\mathbf{d}$ fulfilling (15) and the rank condition of Theil is violated.
Second, if $rk(\mathbf{X}_{\bullet,i}) = K_i$, $\forall i$, then (17) implies the existence
of a nonzero vector $\mathbf{s}$, lying in the linear space spanned by the covariates of equation $i$, for every $i$ with non-zero weight $a_i$ in $\sum_i {a}_i\mathbf{X}_{\bullet,i}\boldsymbol{\beta}_{\bullet,i} = \mathbf{s}$.
Then it is possible to find a nonzero vector $\mathbf{h}_i$ for every $i$ with $a_i \ne 0$ such that
$ \mathbf{X}_{\bullet,i} \mathbf{h}_{i} = \mathbf{s}$.
Finally, by setting $\mathbf{d}_i = a_i \mathbf{h}_i$ for every $i$ with $a_i \ne 0$ and
$\mathbf{d}_i = \mathbf{0}$ for $a_i=0$, again a nonzero $\mathbf{d}$ satisfying (15) is found.
\,\, {$\Box$}
\end{pf}

Remarks.
(a) A Corollary can be found in an early literature and is due to Worswick \& Champernowne (1954): Under homoskedasticity and identical covariates in each equation, Theil's rank condition is never fulfilled, even if there are no collinearities among the $K$ covariates. In this case, the adding-up condition
$\sum_i {a}_i\mathbf{X}_{\bullet,i}\boldsymbol{\beta}_{\bullet,i} = \mathbf{s}$
is automatically fulfilled and hence contains redundant restrictions whenever $T > K$.
(b) Dhrymes \& Schwarz (1987) provide existence conditions for the restricted estimator and show that Theil's rank conditions are either unnecessary or not informative.
Their assumption that all $n$ equations have no common covariate, is a consequence of their construction, assuming that the covariates of each equation are a selection of the $p$ basis vectors.
Our results avoids such unrealistic assumptions and shows that it is sufficient that the covariates lie in a subspace of the linear space spanned by the $p$ basis vectors.

If we additionally wish to impose restrictions (2) and apply the reasoning of the previous subsection we get the TKN-estimator of Theil (1971), Kreijger \& Neudecker (1977),
\begin{equation}
\hat{\boldsymbol{\beta}}_{TKN} = \hat{\boldsymbol{\beta}}_{MLS}
+ \mathbf{C}_+^{-1}\mathbf{R}'\left( \mathbf{R}\mathbf{C}_+^{-1}\mathbf{R}' \right)^{-1} (\mathbf{r-R}\hat{\boldsymbol{\beta}}_{MLS}).
\end{equation}
Let us define a combination of restrictions (2) and (12) such that
\begin{equation}
\mathbf{H} \boldsymbol{\beta} = \mathbf{h},
\end{equation}
where $\mathbf{H}'=(\mathbf{R}' , \mathbf{X}'\mathbf{A})$ and $\mathbf{h}'=(\mathbf{r}' , \mathbf{g}')$.
Then the rank condition which guarantees the existence of the true inverse $(\mathbf{R}\mathbf{C}_+^{-1}\mathbf{R}')^{-1}$ is, in analogy to (4), given by
\begin{equation}
rk(\mathbf{H}) = rk(\mathbf{H},\mathbf{h}),
\end{equation}
which is satisfied whenever (4) holds for a regular design matrix (i.e. $rk(\mathbf{X})=P=K$).

\subsection{Collinear design and singular dispersion}

We treat the multicollinear case $rk(\mathbf{X})=P<K$ together with $rk(\mathbf{\Omega})=M<T$ and consider GLS estimation of the transformed model (10) subject to the restrictions (19). The resulting system of normal equations is given by
\begin{equation} \left( \begin{array}{cc} \mathbf{C}_+ & \mathbf{H}' \\
\mathbf{H} &\mathbf{0} \end{array} \right)
\left( \begin{array}{c} {\hat{\boldsymbol{\beta}}} \\ \boldsymbol{\lambda} \end{array} \right) =
\left( \begin{array}{c} \mathbf{X}' \mathbf{\Omega}^+ \mathbf{y} \\ \mathbf{h} \end{array} \right),
\end{equation}
where $\boldsymbol{\lambda}$ is a vector of Lagrangean multipliers. In contrast to the case of a regular design matrix given in (14), the true inverse of the matrix $\mathbf{C}_+$ in the upper left of the coefficient matrix may not exist for collinear designs. Whenever (4) and (20) hold, system (21) has a solution, as shown above by reparametrization (e.g., Rao, 1965, Haupt \& Oberhofer, 2002).

An alternative approach is to use the base of the null space (e.g., Haupt \& Oberhofer, 2006b) on
$\mathbf{H}$ as the columns of a $(K \times (K-rk(\mathbf{H})))$ matrix $\mathbf{N}$, where without loss of generality we assume $\mathbf{N}'\mathbf{N}=\mathbf{I}_{(K-rk(\mathbf{H}))}$.
Premultiplying the first $K$ rows of (21) by $\mathbf{H}$ and solving for $\boldsymbol{\lambda}$ gives
$\boldsymbol{\lambda} = (\mathbf{HH}')^+ \mathbf{H}
(\mathbf{X}' \mathbf{\Omega}^+ \mathbf{y} - \mathbf{C}_+ {\hat{\boldsymbol{\beta}}})$.
From analogously premultiplying with $\mathbf{N}'$ we get, subject to $\mathbf{HN = 0}$,
\begin{equation}
\mathbf{N}' \mathbf{C}_+ {\hat{\boldsymbol{\beta}}}
= \mathbf{N}' \mathbf{X}' \mathbf{\Omega}^+ \mathbf{y}.
\end{equation}
Since $(\mathbf{H}',\mathbf{N})'$ has full column rank, ${\hat{\boldsymbol{\beta}}}$ can be determined by the equation system constituted by (22) and (19). The unconstrained case is given for an empty
$\mathbf{H}$, where $\mathbf{N = I}_K$. Then (19) represents the usual normal equation system of the Aitken type estimator. If the number of rows in $\mathbf{H}$ exceeds its rank, and (19) and (22) contain more than $K$ equations, some of these are redundant. See Markiewicz \& Puntanen (2015) for an extensive treatment and literature survey.

\begin{thm}
For the regression model defined in (1), equations (19) and (22) have at most one solution,
whenever (3), (4), (7), and (20) hold.
\end{thm}

\begin{pf}
We prove Theorem 2 by contradiction.
Let $\Delta \boldsymbol{\beta}$  denote the difference between two different solutions.
Then the equations
$\mathbf{N}' \mathbf{C}_+ \Delta \boldsymbol{\beta} = \mathbf{0}$ and
$\mathbf{H} \Delta \boldsymbol{\beta} = \mathbf{0}$
must be fulfilled, respectively.
Due to $\mathbf{HN = 0}$ follows the existence of a vector $\mathbf{c}$, which fulfills
$\Delta \boldsymbol{\beta} = \mathbf{Nc}$. Thus
$\mathbf{N}' \mathbf{C}_+ \mathbf{Nc} = \mathbf{0}$.
According to Haupt \& Oberhofer (2002, Lemma 1), the matrix
$\mathbf{S} = \mathbf{N}'\mathbf{C}_+ \mathbf{N}$ of order $K -rk(\mathbf{H})$ is invertible.
Thus, $\Delta \boldsymbol{\beta}$ has to be equal to zero.
\,\, {$\Box$}
\end{pf}

\begin{thm} Under the assumptions of Theorem 2 and an arbitrary vector $\boldsymbol{\beta}^*$
fulfilling $\mathbf{H}\boldsymbol{\beta}^* = \mathbf{h}$, equations (19) and (22) have a unique solution
\begin{equation}
{\hat{\boldsymbol{\beta}}} = \mathbf{NS}^{-1}\mathbf{N}'\mathbf{X}' \mathbf{\Omega}^+ \mathbf{y} +
(\mathbf{I}_K - \mathbf{NS}^{-1} \mathbf{N}'\mathbf{C}_+) \boldsymbol{\beta}^*,
\end{equation}
where  $E(\hat{\boldsymbol{\beta}})=\boldsymbol{\beta}$ and $Var(\hat{\boldsymbol{\beta}})=\sigma^2 \mathbf{NS}^{-1}\mathbf{N}'$.
\end{thm}

\begin{pf}
Due to Theorem 2 it is sufficient to show that (23) is a solution to (22) and satisfies (21).
This follows by substituting the estimator (23) into (19) and the definition of $\mathbf{S}$.
From (23) and (19) follows
$\mathbf{H}\hat{\boldsymbol{\beta}} = \mathbf{H}\boldsymbol{\beta}^* = \mathbf{h}$.
\,\, {$\Box$}
\end{pf}

Remark.
Note that $\hat{\boldsymbol{\beta}}$ is unique but, due to the adding-up restriction, does not have a unique representation. If the parameter space is constrained, an optimal estimator is always affine
$\hat{\boldsymbol{\beta}} = \mathbf{Ly + l}$.
Representing $\hat{\boldsymbol{\beta}}$ as a function of $\mathbf{y}$ enables an arbitrary interchange between the homogenous and the particular part of the estimator, leading to a class of linear representations. The following general result is due to Haupt \& Oberhofer (2006b).

\begin{thm}
Under the assumptions of Theorem 2,
${\hat{\boldsymbol{\beta}}}$ in (23) is the best affine unbiased estimator of $\boldsymbol{\beta}$ in (1). The class of all linear
representations of ${\hat{\boldsymbol{\beta}}}$ is given by
\begin{equation}
{\hat{\boldsymbol{\beta}}} =
(\mathbf{NS}^{-1}\mathbf{N}'\mathbf{X}' \mathbf{\Omega}^+ + \mathbf{GA}')\mathbf{y} +
(\mathbf{I}_K - \mathbf{NS}^{-1} \mathbf{N}'\mathbf{C}_+) \boldsymbol{\beta}^* - \mathbf{Gg},
\end{equation}
due to $\mathbf{A}' \mathbf{y}=\mathbf{g}$ with probability one, where $\mathbf{G}$ is an arbitrary matrix.
\end{thm}

\begin{pf}
See Haupt \& Oberhofer (2006b, Theorems 1 and 2). \,\, {$\Box$}
\end{pf}

\section{FE panel regression and the within transformation}

A special case of model class (SUR) are so-called fixed effects (FE) panel regression models. They are based on the assumption that the parameters do not vary over equations $i$.
The one-way FE model includes a covariate meant to capture the unobserved equation-specific heterogeneity, which is fixed over time (e.g., Kiefer, 1980, Im et al., 1999),
\begin{equation*}
\mu_{t,i}=\sum_{k=1}^K \beta_{k} x_{t,k,i} + \gamma_i,
\quad 1 \le i \le n, \quad 1 \le t \le m.
\tag{FE}
\end{equation*}
A more general approach is the two-way FE model
$\mu_{t,i}=\sum_{k=1}^K \beta_{k} x_{t,k,i} + \gamma_i + \delta_t$, $1 \le i \le n$, $1 \le t \le m$,
where $\delta_t$ is included to capture unobserved time-specific effects common to all equations.
For every $n$, model (FE) can be written as
\[
\mathbf{y}_{\bullet,i} = \mathbf{X}_{\bullet,i} \boldsymbol{\beta} + \mathbf{e}_m \gamma_i  + \mathbf{u}_{\bullet,i},
\quad 1 \le i \le n,
\]
where $\mathbf{e}_{m}$ is the $(m \times 1)$ vector of ones,
$\mathbf{y}_{\bullet,i}$ is the $(m \times 1)$ vector of response variables in equation $i$,
$\mathbf{X}_{\bullet,i}$ is an $(m \times K)$ covariate matrix, $\boldsymbol{\beta}$
contains the $K$ parameters, and $\mathbf{u}_{\bullet,i}$ is an $(m \times 1)$ error vector.
In compact notation and under the assumptions stated after equation (1),
\begin{equation}
\mathbf{y} = \mathbf{X} \boldsymbol{\beta} + \mathbf{Z} \boldsymbol{\gamma} + \mathbf{u},
\label{eq:FE}
\end{equation}
where all dimensions equal those in equation (1), $\boldsymbol{\gamma}=({\gamma}_1, \ldots, {\gamma}_n)'$,
and $\mathbf{Z}=\mathbf{I}_n \otimes \mathbf{e}_m$ is the $(T \times n)$ matrix of equation specific dummy variables.

Due to the so-called incidental parameter problem, the literature discussing equation (25) is focussed on partial estimation of the subset of parameters $\boldsymbol{\beta}$. Under a suitable rank condition on $(\mathbf{X},\mathbf{Z})$, as stated in Section 2, the OLS estimator of $\boldsymbol{\beta}$ in equation (25) is linear unbiased and given by
\begin{equation*}
\hat{\boldsymbol{\beta}}_{OLS}
=(\mathbf{X}'\mathbf{M}\mathbf{X})^{-1}
\mathbf{X}'\mathbf{M}\mathbf{y},
\end{equation*}
where $\mathbf{M}$ is the projection matrix
$\mathbf{M}
\equiv \mathbf{I}_T - \mathbf{Z}(\mathbf{Z}'\mathbf{Z})^{-1}\mathbf{Z}'
= \mathbf{I}_n \otimes \mathbf{M}_m$,
$\mathbf{M}_m \equiv \mathbf{I}_m - \mathbf{e}_m\mathbf{e}_m'/m$
is the so-called centering matrix, and thus a typical element of the $(T \times 1)$ vector $\mathbf{M}\mathbf{y}$ is equal to $y_{t,i}-\bar{y}_{i}$.

If the error variables are assumed to be heteroskedastic over equations and time,
we have
$E(\mathbf{u}_{\bullet,i}\mathbf{u}_{\bullet,i}') \equiv \mathbf{\Sigma}_{i,i}$, $1 \le i \le n$,
and $\mathbf{\Sigma}_{i,j}=\mathbf{0}$, for $i \neq j$,
where $\mathbf{\Sigma}_{1,1}, \ldots, \mathbf{\Sigma}_{n,n}$ are known symmetric, nonnegative definite $(m \times m)$ matrices with elements
$Cov(u_{s,i},u_{t,i}) \equiv
\sigma_{s,t,i,i}$, $1 \le s,t \le m, \,\,1 \le i \le n$.
Then,
$E(\mathbf{{u}} \mathbf{{u}})' \equiv \sigma^2 \boldsymbol{\Omega} = \diag{\mathbf{\Sigma}_{1,1}, \ldots, \mathbf{\Sigma}_{n,n}}$,
with $rk(\boldsymbol{\Omega})=T$.
For ease of exposition we consider a system dispersion matrix with Kronecker structure
$E(\mathbf{{u}} \mathbf{{u}}') = \mathbf{I}_n \otimes \mathbf{\Sigma}$ (e.g., Kiefer, 1980),
under the assumption of homoskedasticity over equations but arbitrary intertemporal covariance (i.e. $\sigma_{s,t,1,1} = \ldots = \sigma_{s,t,n,n} = \sigma_{s,t}$).
All of the results given below hold for general block-diagonal dispersion.

The GLS estimator of $\boldsymbol{\beta}$ in (25) is best linear unbiased and given by
\begin{equation}
\hat{\boldsymbol{\beta}}_{GLS}
=(\mathbf{X}'\mathbf{P}\mathbf{X})^{-1}
\mathbf{X}'\mathbf{P}\mathbf{y},
\end{equation}
where
$\mathbf{P} \equiv (\mathbf{I}_n \otimes \mathbf{\Sigma}^{-1})(\mathbf{I}_T - \mathbf{Q})$
and
$\mathbf{Q} \equiv \mathbf{Z}(\mathbf{Z}'(\mathbf{I}_n \otimes \mathbf{\Sigma}^{-1})\mathbf{Z})^{-1}\mathbf{Z}' (\mathbf{I}_n \otimes \mathbf{\Sigma}^{-1})$
is idempotent.
Obvious consequences of these definitions are $\mathbf{Q}\mathbf{Z} = \mathbf{Z}$, $\mathbf{P}\mathbf{Z} = \mathbf{0}$,
$\mathbf{P}\mathbf{M} = \mathbf{P}$, $\mathbf{M}\mathbf{P}\mathbf{M} = \mathbf{P}$, and
$\mathbf{P}(\mathbf{I}_n \otimes \mathbf{\Sigma})\mathbf{P} = \mathbf{P}$.

A generic example for regression systems with regular design matrix and singular dispersion matrix discussed in Section 2.2 results from applying the within transformation (e.g., Kiefer, 1980, Im et al., 1999, Qian \& Schmidt, 2003, or Tian \& Jiang, 2016).
It is well known that $\boldsymbol{\beta}$ can be estimated (consistently for large $n$ and fixed $m$) by applying OLS to the transformed model
\begin{equation}
\mathbf{M}\mathbf{y}=\mathbf{M}\mathbf{X} \boldsymbol{\beta} + \mathbf{M}\mathbf{u}.
\end{equation}
Due to $rk(\mathbf{M}_m)=m-1$, the error term $\mathbf{M}\mathbf{u}$ has a singular dispersion matrix
$E(\mathbf{M}\mathbf{uu}'\mathbf{M})
= E(\mathbf{I}_n \otimes \mathbf{M}_m\mathbf{\Sigma}\mathbf{M}_m)$,
where $rk(\mathbf{M}_m\mathbf{\Sigma}\mathbf{M}_m)=T-1$.

Instead of dealing with the problem of the singular dispersion matrix when applying GLS to equation (27), it is well known that dropping one time period leads to algebraically identical results (e.g., Kiefer, 1980). In the following we provide a rigorous proof of this assertion using the results from the previous section.
We proceed as discussed in Section 2.2 and consider the spectral decomposition
\[
\mathbf{I}_n \otimes \mathbf{M}_m\mathbf{\Sigma}\mathbf{M}_m
= (\mathbf{a},\mathbf{F}) \diag{{0}, \mathbf{\Lambda}} (\mathbf{a}',\mathbf{F}')'
= \mathbf{F}\mathbf{\Lambda}\mathbf{F}',
\]
where, without loss of generality, we choose the $(T \times 1)$ vector $\mathbf{a}$ to be the first column of the ($T \times m$) matrix
$\mathbf{A}
= (\mathbf{a}_{(1)}, \ldots, \mathbf{a}_{(t)}, \ldots, \mathbf{a}_{(m)})
= \mathbf{I}_n \otimes \mathbf{e}_m/\sqrt{m}$.
Premultiplying equation (27) with $\mathbf{a}_{(1)}'$ is equal to eliminating all observations corresponding to $t=1$.
Then, $E(\mathbf{a}'\mathbf{M}\mathbf{uu}'\mathbf{M}\mathbf{a}) = \sigma^2\mathbf{a}'\mathbf{M}\boldsymbol{\Omega}\mathbf{M}\mathbf{a} = 0$, and
according to (11), we have $\mathbf{a}'\mathbf{M}\mathbf{y}=\mathbf{a}'\mathbf{M}\mathbf{X} \boldsymbol{\beta}$ with probability one.

From Theorem 1 of Haupt \& Oberhofer (2006a) follows that we obtain the same results when we delete observation $t=1$ with $\mathbf{a}_{(1)} \neq \mathbf{0}$, or solve the normal equations
\[
\mathbf{X}'\mathbf{M}(\mathbf{I}_n \otimes \mathbf{M}_m\mathbf{\Sigma}\mathbf{M}_m)^{+}\mathbf{M}\mathbf{X}\boldsymbol{\beta}
= \mathbf{X}'\mathbf{M}(\mathbf{I}_n \otimes \mathbf{M}_m\mathbf{\Sigma}\mathbf{M}_m)^{+}\mathbf{y},
\]
using the Moore-Penrose inverse
$(\mathbf{I}_n \otimes \mathbf{M}_m\mathbf{\Sigma}\mathbf{M}_m)^{+} = \mathbf{F}\mathbf{\Lambda}^{-1}\mathbf{F}'$.
Then, in analogy to (14) and from Theorem 1, follows that the true inverse of $\mathbf{X}'\mathbf{M}(\mathbf{I}_n \otimes  \mathbf{M}_m\mathbf{\Sigma}\mathbf{M}_m)^{+}\mathbf{M}\mathbf{X}$ exists.
Hence, the GLS estimator of $\boldsymbol{\beta}$ in equation (27) exists and is equal to the Moore-Penrose inverse least squares estimator (discussed in Section 2.2)
\begin{equation}
\hat{\boldsymbol{\beta}}_{MLS} = \mathbf{C}_+^{-1} \mathbf{X}'\mathbf{M}(\mathbf{I}_n \otimes \mathbf{M}_m\mathbf{\Sigma}\mathbf{M}_m)^{+}\mathbf{y},
\end{equation}
where
$\mathbf{C}_+^{-1} = (\mathbf{X}'
\mathbf{M}(\mathbf{I}_n \otimes \mathbf{M}_m\mathbf{\Sigma}\mathbf{M}_m)^{+}\mathbf{M}
\mathbf{X})^{-1}
= (\mathbf{X}'\mathbf{M} \mathbf{F} \mathbf{\Lambda}^{-1} \mathbf{F}'\mathbf{M}\mathbf{X})^{-1}$.

\begin{thm}
Let $\hat{\boldsymbol{\beta}}_{GLS}$ as defined in equation (26) and $\hat{\boldsymbol{\beta}}_{MLS}$ as defined in equation (28). Then
$\hat{\boldsymbol{\beta}}_{MLS} = \hat{\boldsymbol{\beta}}_{GLS}$.
\end{thm}

\begin{pf}
We have to show
$\mathbf{P}=\mathbf{M}(\mathbf{I}_n \otimes \mathbf{M}_m\mathbf{\Sigma}\mathbf{M}_m)^{+}\mathbf{M}$,
which implies
$\mathbf{P}=(\mathbf{I}_n \otimes \mathbf{M}_m\mathbf{\Sigma}\mathbf{M}_m)^{+}$.
Using the spectral decomposition it remains to verify
$\mathbf{P}\mathbf{F}\mathbf{\Lambda}\mathbf{F}' = \mathbf{F}\mathbf{\Lambda}^{-1}\mathbf{F}'\mathbf{F}\mathbf{\Lambda}\mathbf{F}' = \mathbf{I}_T$.
Using the results following equation (26), we have
\[
\mathbf{P}\mathbf{F}\mathbf{\Lambda}\mathbf{F}'
=\mathbf{P}(\mathbf{I}_n \otimes \mathbf{M}_m\mathbf{\Sigma}\mathbf{M}_m)
=\mathbf{P}\mathbf{M}(\mathbf{I}_n \otimes \mathbf{\Sigma})\mathbf{M} = \mathbf{I}_T.
\]
This implies $\mathbf{P}(\mathbf{I}_n \otimes \mathbf{\Sigma})\mathbf{P}=\mathbf{P}$.
\,\, {$\Box$}
\end{pf}

Remarks.
(a) The previous results (and Theorem 1 of Haupt and Oberhofer, 2006a) can be applied to the more general panel data model with two-way fixed-effects
$\mu_{t,i}=\sum_{k=1}^K \beta_{k} x_{t,k,i} + \gamma_i + \delta_t$, $1 \le i \le n$, $1 \le t \le m$.
(b) Qian \& Schmidt (2003, Theorem 1) show that the GLS estimator of $\boldsymbol{\beta}$ in equation (27) is not just an alternative to the GLS estimator of $\boldsymbol{\beta}$ in equation (25), but that these estimators are, in fact, equal. Their result is embedded by many previous results on general Gauss-Markoff theory (see, e.g., Tian and Zhang, 2011). It is of particular interest for the estimation of the Kronecker-type dispersion matrix when $n$ goes to infinity for fixed $m$.
(c) As stated above, Theorem 5 holds for any block-diagonal dispersion structure.
Deriving Theorem 5 under the collinearity assumptions of Section 2.3 is left to the reader.

\end{document}